\documentclass[11pt,dvips,epsfig,reqno,graphics]{amsart}
\usepackage{a4wide, latexsym,amssymb,amsthm,amsmath,amsfonts,url}

\begin{document}

\title{A new approach to gradient Ricci solitons and generalizations}
\author{M. Crasmareanu}
\date{July 11, 2017}
\dedicatory{Dedicated to Academician Radu Miron on the occasion of his 90'th birthday}

\begin{abstract}
This short note concerns with two inequalities in the geo\-me\-try of gradient Ricci solitons $(g, f, \lambda )$ on a smooth manifold $M$. These inequalities provide some relationships
between the curvature of the Riemannian metric $g$ and the behavior of the scalar field $f$ through two second order equations satisfied by the scalar $\lambda $. We propose several
ge\-ne\-ra\-li\-za\-tions of Ricci solitons to the setting of manifolds endowed with linear connections, not necessary of metric type.
\end{abstract}



\maketitle

\medskip

Let $(M^n, g)$ be a $n$-dimensional Riemannian manifold endowed with a smooth function $f\in C^{\infty }(M)$. The scalar field $f$ yields {\it the Hessian endomorphism}:
$$
h_f:\mathcal{X}(M)\rightarrow \mathcal{X}(M), \quad h_f(X)=\nabla _X\nabla f \eqno(1)
$$
where $\nabla $ is the Levi-Civita connection of $f$. Then we know the symmetry of {\it the Hessian tensor field} of $f$:
$$
H_f(X, Y):=g(h_f(X), Y), \quad H_f(X, Y)=H_f(Y, X) \eqno(2)
$$
It follows the existence of a $g$-orthonormal frame $E=\{E_i\}_{i=1,...,n}\in \mathcal{X}(M)$ and the existence of the eigenvalues $\lambda =\{\lambda _i\}_{i=1,...,n}\in C^{\infty
}(M)$:
$$
h_f(E_i)=\lambda _iE_i. \eqno(3)
$$
Hence we express all the geometric objects related to $f$ in terms of the pair $(E, \lambda )$:
$$
\nabla f=\sum _{i=1}^nE_i(f)E_i, \quad \|\nabla f\|_g^2=\sum _{i=1}^n[E_i(f)]^2, \quad h_f(X)=\nabla _X\nabla f=\sum _{i=1}^n(\lambda _iX^i)E_i, \eqno(4)
$$
for $X=X^iE_i$. Also the Hessian and {\it the Laplacian} of $f$ are:
$$
H_f(X, Y)=\sum _{i=1}^n\lambda _i(X^iY^i), \quad \Delta f:=Tr_g H_f=\sum _{i=1}^n\lambda _i. \eqno(5)
$$
Suppose now that the triple $(g, f, \lambda \in \mathbb{R})$ is a {\it gradient Ricci soliton} on $M$, \cite[p. 76]{c:g}:
$$
H_f+Ric+\lambda g=0 \eqno(6)
$$
where $Ric$ is the Ricci tensor field of $g$. By considering the Ricci endomorphism $Q\in \mathcal{T}^1_1(M)$ provided by:
$$
Ric(X, Y)=g(QX, Y) \eqno(7)
$$
we can express $(6)$ as:
$$
h_f+Q+\lambda I=0 \eqno(8)
$$
with $I$ being the Kronecker endomorphism. From $(4)$ we get that $Q$ is also of diagonal form with respect to the frame $E$:
$$
Q(X)=-\sum_{i=1}^n(\lambda _i+\lambda )X^iE_i, \quad \|Q\|_g^2=\sum _{i=1}^n(\lambda _i+\lambda )^2. \eqno(9)
$$
By developing the second formula above we derive:
$$
\|Ric\|_g^2=\sum _{i=1}^n\lambda _i^2+2\lambda \sum _{i=1}^n\lambda _i+n\lambda ^2=\|H_f\|_g^2+2\lambda \Delta f+n\lambda ^2. \eqno(10)
$$
Hence the scalar $\lambda $ is a solution of the second order equation:
$$
n\lambda ^2+2\Delta f\cdot \lambda +(\|H_f\|_g^2-\|Ric \|_g^2)=0 \eqno(11)
$$
which means the positivity:
$$
0\leq \Delta ^{\prime}:=(\Delta f)^2-n(\|H_f\|_g^2-\|Ric \|_g^2). \eqno(12)
$$
It follows a boundary of the geometry of $g$ in terms of $f$:
$$
\|Ric\|_g^2\geq \|H_f\|_g^2-\frac{1}{n}(\Delta f)^2. \eqno(13)
$$
An "exotic" consequence is provided by the case of strict inequality in $(12)$; more precisely, it follows that the data $(g, f, \lambda )$ is doubled by the triple $(g, f,
-\frac{2\Delta f}{n}-\lambda )$.

\medskip

{\bf Examples} 1) (Gaussian soliton) We have $(M=\mathbb{R}^n, g_{can})$ and $f(x)=-\frac{\lambda }{2}\|x\|^2$. It results: $h_f=-\lambda I_n$ and $\Delta f=-n\lambda $. Since
$\|H_f\|^2=n\lambda ^2$ the left hand side of $(11)$ is:
$$
n\lambda ^2+2\Delta f\cdot \lambda +(\|H_f\|_g^2-\|Ric \|_g^2)=n\lambda ^2+2(-n\lambda )\lambda +(n\lambda ^2-0)
$$
which is exactly zero. Also: $ \Delta ^{\prime}=(n\lambda )^2-n(n\lambda ^2-0)=0$ which means the uniqueness of $\lambda $ and the equality case in $(13)$:
$0=n\lambda ^2-\frac{(n\lambda )^2}{n}$. \\
2) (Einstein manifold) Let $(M, g)$ be an Einstein manifold: $Ric =-\lambda g$. A function $f$ with vanishing Hessian is called {\it Killing potential} in \cite{c:m1} since its gradient
is a Killing vector field; in \cite[p. 283]{p:p} such a function is called {\it linear}. Hence: $\Delta f=\|H_f\|^2=0$ and $\|Ric\|^2=n\lambda ^2$ which yields the following value of the
left hand side of $(11)$:
$$
n\lambda ^2+2\Delta f\cdot \lambda +(\|H_f\|_g^2-\|Ric \|_g^2)=n\lambda ^2+2(\cdot 0)\lambda +(0-n\lambda ^2)
$$
which is exactly zero. Also: $ \Delta ^{\prime}=0^2-n(0-n\lambda ^2)=n^2\lambda ^2\geq 0$ which means that $(13)$ is the inequality $n\lambda ^2\geq 0-\frac{0}{n}$ and the uniqueness of
$\lambda $ gives a steady soliton, equivalently $g$ is Ricci-flat. We consider an in\-te\-res\-ting open problem to find the linear functions of a steady
soliton and of a Ricci-flat metric. \\
3) A generalization of the first example is provided on a Ricci-flat manifold by a smooth function $f$ satisfying a generalization of Hessian structures:
$$
H_f=-\lambda g. \eqno(14)
$$
Then: $\Delta f=-n\lambda $ and $\|H_f\|^2=n\lambda ^2$ exactly as for the Gaussian soliton. Using Lemma 4.1. of \cite[p. 1540]{c:c} it results form $(14)$ that $\nabla f$ is a
particular {\it concircular vector field}: $h_f=-\lambda I$; hence $\lambda _1=...=\lambda _n=-\lambda $. If $\nabla f$ is without zeros it follows from Theorem 3.1. of \cite[p.
1539]{c:c} that $(M, g)$ is locally a warped product with a $1$-dimensional basis: $(M, g)=(I\times _{\varphi }(F^{n-1}, g_F))$. In fact: $\nabla f=\varphi (s)\frac{\partial }{\partial
s}$ with $\varphi ^{\prime }(s)=-\lambda $ which means a affine warping function: $\varphi (s)=-\lambda s+C$. \\
4) (Hamilton's cigar) The famous Hamilton's soliton is the steady soliton provided by the complete Riemannian geometry $(\mathbb{R}^2, g=\frac{1}{1+x^2+y^2}g_{can})$ and the potential
function $f(x, y)=-\ln (1+x^2+y^2)$. The only non-zero components of the Hessian are $(H_f)_{11}=(H_f)_{22}=\frac{-2}{(1+x^2+y^2)^2}$ which yields the norm
$\|H_f\|^2=\frac{8}{(1+x^2+y^2)^2}$ and the Laplacian $\Delta f=\frac{-4}{1+x^2+y^2}$. The Gaussian curvature of $g$ is $K(x, y)=\frac{2}{1+x^2+y^2}$ and $\|Ric\|^2=2K^2$. In conclusion,
the inequality $(13)$ becomes strictly: $\frac{8}{(1+x^2+y^2)^2}>\frac{8}{(1+x^2+y^2)^2}-\frac{16}{2(1+x^2+y^2)^2}=0$. \\
5) (Cylinder shrinking soliton) Consider the Riemannian product $(M, g)=(\mathbb{S}^{n-1}\times \mathbb{R}, 2(n-2)g_S+dt^2)$ of the shrinking unit sphere for $n\geq 3$ with a line. We
have:
$$
Ric_g=(n-2)g_S=\frac{1}{2}(g-dt^2), \quad \|Ric\|^2=\frac{1}{4}. \eqno(15)
$$
For $f(x, t)=\frac{t^2}{4}$ we obtain a shrinking gradient Ricci soliton with $\lambda =-\frac{1}{2}$ from $H_f=\frac{1}{2}dt^2$; also: $\Delta f=\frac{1}{2}$. Then (13) becomes the
strict inequality: $\frac{1}{4}>\frac{1}{4}-\frac{1}{4n}$.  \quad $\Box $

\medskip

A new second order equation, similar to $(11)$, follows from a well-known formula from the theory of gradient Ricci solitons, \cite[p. 79]{c:g}:
$$
\Delta f+R+n\lambda =0 \eqno(16)
$$
obtained by tracing $(6)$; here $R$ is the scalar curvature of $g$. Hence the companion equation of $(11)$ is:
$$
n\lambda ^2+2R\lambda +(\|Ric\|^2_g-\|H_f\|^2_g)=0. \eqno(17)
$$
The new inequality is then:
$$
0\leq \Delta ^{\prime }:=R^2-n(\|Ric\|^2_g-\|H_f\|^2_g) \eqno(18)
$$
and it results a boundary of the behavior of $f$ in terms of geometry of $g$:
$$
\|H_f\|^2_g\geq \|Ric\|^2_g-\frac{R^2}{n}. \eqno(19)
$$
We remark that $(13)$ and $(19)$ can be unified in the double inequality:
$$
\|H_f\|^2_g-\frac{1}{n}(\Delta f)^2\leq \|Ric\|^2_g\leq \|H_f\|^2_g+\frac{R^2}{n} \eqno(20)
$$
and the simultaneous equalities hold if and only if: $H_f=\pm Ric$ with $\Delta f=R=0$; hence $f$ is a harmonic map.

\smallskip

{\bf Examples revisited} 1) (Gaussian soliton) The inequality $(19)$ becomes: $n\lambda ^2\geq 0$. \\
2) (Einstein manifold) $(19)$ becomes an equality: $0=n\lambda ^2-\frac{(n\lambda )^2}{n}$. \\
3) (Hamilton soliton) $(19)$ is the strict inequality: $\frac{8}{(1+x^2+y^2)^2}>\frac{2\cdot 4}{(1+x^2+y^2)^2}-\frac{4^2}{2(1+x^2+y^2)^2}=0$. \\
4) (Cylinder shrinking soliton) Since $R=\frac{n-2}{2}$ we get the following form of $(19)$: $\frac{1}{4}>\frac{1}{4}-\frac{(n-2)^2}{4n^2}$. \quad $\Box $

\medskip

In the second part of this note we connect the above considerations with a study of Academician Radu Miron who was the Adviser of the author's PhD Thesis in February 1999. More
precisely, let $\tilde{\nabla }$ a linear connection different to the Levi-Civita connection of $g$. Let also the $1$-form $\eta =df$ be the $g$-dual of $\xi =\nabla f$. The data $(M, g,
\tilde{\nabla }, \eta )$ is called {\it Weyl space} in \cite{m:r} if $g$ is $\tilde{\nabla }$-{\it recurrent} with the factor $\eta $:
$$
\tilde{\nabla }g=\eta \otimes g. \eqno(21)
$$
Hence, we arrive at a generalization of Ricci solitons in the framework of manifolds endowed with a linear connection:

\medskip

{\bf Definition 1} Let $(M, \tilde{\nabla })$ be given. A triple $(g, \xi \in \mathcal{X}(M), \lambda )$ is a $\tilde{\nabla }$-{\it Ricci soliton} if:
$$
\tilde{\nabla }\xi +Q+\lambda I=0. \eqno(22)
$$
Let $\eta =\xi ^{\flat}$ be the $g$-dual form of $\xi $. The data $(g, \xi , \lambda , \mu \in \mathbb{R})$ is a $(\tilde{\nabla }, \eta )$-{\it Ricci soliton} if:
$$
\tilde{\nabla }\xi +Q+\lambda I+\mu \eta \otimes \xi =0. \eqno(23)
$$
More generally, let $(M, \tilde{\nabla }, F\in \mathcal{T}^1_1(M))$ be given. The pair $(\xi , \lambda )$ will be a $(\tilde{\nabla }, F)$-{\it soliton} if:
$$
\tilde{\nabla }\xi +F+\lambda I=0 \eqno(24)
$$
and the triple $(\xi , \lambda , \mu )$ is a $(\tilde{\nabla }, F, \eta )$-{\it soliton} if:
$$
\tilde{\nabla }\xi +F+\lambda I+\mu \eta \otimes \xi =0. \eqno(25)
$$

\smallskip

On this way we propose a study of classes of solitons, maybe more adapted to Hermitian/K\"ahler geometry ($F^2=-I$) and para-Hermitian/para-K\"ahler geometry ($F^2=I$) by using some
linear connections adapted to these settings like the Chern and Bismut complex connections, \cite{c:m2}.

\medskip

As first example, we consider the Vaisman geometry following \cite{b:c}. Let $(M^{2n}, J, g)$ be a complex $n$-dimensional Hermitian manifold and $\Omega $ its fundamental $2$-form given
by $\Omega (X, Y)=g(X, JY)$ for any vector fields $X, Y\in \Gamma (TM)$. Recall from \cite[p. 1]{d:ob} that $(M, J, g, \Omega )$ is a {\it locally conformal K\"ahler manifold} (l.c.K) if
there exists a closed $1$-form $\omega \in \Gamma (T^0_1(M))$ such that: $d\Omega =\omega \wedge \Omega $. In particular, $M$ is called {\it strongly non-K\"ahler} if $\omega $ is
without singularities i.e. $\omega \neq 0$ everywhere; hence we consider $2c=\|\omega \|$ and $u=\omega /2c$ the corresponding $1$-form. Since $\omega $ is called {\it the Lee form} of
$M$ the vector field $U=u^{\sharp}$ will be called {\it the Lee vector field}. Consider also the unit vector field $V=JU$, {\it the anti-Lee vector field}, as well as its dual form
$v=V^{\flat}$, so: $u(V)=v(U)=0, v=-u\circ J, u=v\circ J$.

\medskip

Our setting is provided by the particular case of strongly non-K\"ahler l.c.K. manifolds, called {\it Vaisman manifolds}, and given by the parallelism of $\omega $ with respect to the
Levi-Civita connection $\nabla $ of $g$. Hence $c$ is a positive constant and the Lemma 2 of \cite{n:p} gives the covariant derivative of $V$ with respect to any $X\in \Gamma (TM)$:
$$
\nabla _XV=c[u(X)V-v(X)U-JX]. \eqno(26)
$$
It follows a class of general solitons provided by:

\medskip

{\bf Proposition 1} {\it Let $(M, J, g, c=1)$ be a Vaisman manifold and the linear connection:
$$
\tilde{\nabla }:=\nabla -[u\otimes I+v\otimes J]. \eqno(27)
$$
Then $(V, 0\in \mathbb{R})$ is a $(\tilde{\nabla }, J)$-soliton. Moreover, $g$ is recurrent with respect to $\tilde {\nabla }$ with the factor $2u$:
$$
\tilde {\nabla }g=2u\otimes g \eqno(28)
$$
but $\tilde{\nabla }$ is not the corresponding Weyl connection since is not torsion-free}:
$$
\tilde{T}=I\otimes u+J\otimes v-u\otimes I-v\otimes J. \eqno(29)
$$

\smallskip

As second example we consider $\tilde{\nabla }$ as being exactly the Weyl connection $(21)$ of the pair $(g, \eta)$. Is well-known its expression:
$$
\tilde{\nabla }=\nabla -\frac{1}{2}\eta \otimes I-\frac{1}{2}I\otimes \eta +\frac{1}{2}g\otimes \xi \eqno(30)
$$
and we derive:

\medskip

{\bf Proposition 2} {\it Let $(M, g, \eta , \tilde{\nabla })$ be a Weyl geometry with $\eta $ of constant norm and endowed with an endomorphism $F$. Then $(\xi =\eta ^{\sharp}, \lambda ,
\mu )$ is a $(\tilde{\nabla }, F, \eta )$-soliton if and only if $(\xi, \lambda -\frac{\|\xi \|^2_g}{2}, \mu +\frac{1}{2})$ is a $(\nabla , F, \eta )$-soliton with $\nabla $ the usual
Levi-Civita connection of} $g$.

\medskip

As third example, we consider the vector field $\xi $ as  being \textit{torse-forming} on the Riemannian manifold $(M, g)$:
$$
\nabla \xi =fI+\gamma \otimes \xi \eqno(31)
$$
for a smooth function $f\in C^{\infty}(M)$ and a $1$-form $\gamma \in \Omega ^1(M)$. Note that torse-forming vector fields appear in many areas of differential geometry and physics as is
pointed out in \cite{a:im} and are natural generalizations of concircular vector fields. We get immediately:

\medskip

{\bf Proposition 3} {\it Suppose that $\xi $ is a special torse-forming vector field having $f$ a constant function and $\gamma =\eta =\xi ^{\flat}$. Then $(g, \xi , \lambda )$ is a
Ricci soliton if and only if $g$ is an eta-Einstein metric}:
$$
Q=-(\lambda +f)I-\eta \otimes \xi . \eqno(32)
$$

\medskip

As fourth example, let $(M, g, \nabla )$ be a hypersuface of the Riemannian manifold $(\tilde{M}^{n+1}, \tilde{g}, \tilde{\nabla })$ and $A$ its shape operator.

\medskip

{\bf Definition 2} The pair $(\xi \in \mathcal{X}(M), \lambda )$ is {\it a shape soliton on} $(M, g)$ if:
$$
\nabla \xi +A+\lambda I=0. \eqno(33)
$$

\smallskip

{\bf Remark 1} For example, if $M$ is eta-umbilical i.e. $A$ is of eta-type which means that it is has two eigenvalues: $A=\sigma I+(\rho -\sigma )\eta \otimes \xi $ then the above
condition yields that $\xi $ is torse-forming:
$$
\nabla \xi =-(\lambda +\sigma )I+(\sigma -\rho )\eta \otimes \xi . \eqno(34)
$$
Hopf hypersurfaces of eta-umbilical type are studied in \cite[p. 60]{c:c}. Also, the CR submanifolds of maximal CR dimension in complex projective space have eta-type shape operators as
is pointed out in \cite[p. 190]{c:p}. \quad $\Box $

\medskip

We go further with a generalization to the setting of statistical structures of \cite{c:u} provided by data $(M, g, \tilde{\nabla }, \tilde{\nabla }^{\ast })$ where $\tilde{\nabla }$,
$\tilde{\nabla }^{\ast}$ is a pair of torsion-free {\it dual connections} on $(M, g)$:
$$
Z(g(X, Y))=g(\tilde{\nabla }_ZX, Y)+g(X, \tilde{\nabla }^{\ast }_ZY) \eqno(35)
$$
for any vector fields $X$, $Y$, $Z$. We introduce:

\medskip

{\bf Definition 3} i) The statistical manifold $(M, g, \tilde{\nabla }, \tilde{\nabla }^{\ast })$ is called {\it Ricci-symmetric} if the Ricci tensor field $\tilde{Q}$ of $\tilde{\nabla
}$ (equivalently, $\tilde{Q}^{\ast}$ of $\tilde{\nabla }^{\ast}$ by Corollary 9.5.3 of \cite[p. 267]{c:u}) is symmetric. \\
ii) The pair $(\xi , \lambda )$ is a {\it statistical soliton} for the Ricci-symmetric statistical manifold $(M, g, \tilde{\nabla }, \tilde{\nabla }^{\ast })$ if the triple $(g, \xi ,
\lambda )$ is both $\tilde{\nabla }$-Ricci and $\tilde{\nabla }^{\ast}$-Ricci soliton.

\medskip

{\bf Remark 2} Since the Levi-Civita connection of $g$ is the arithmetic mean of the pair $(\tilde{\nabla }, \tilde{\nabla }^{\ast })$:
$$
2\nabla =\tilde{\nabla }+\tilde{\nabla }^{\ast } \eqno(36)
$$
it follows that a statistical soliton is an usual Ricci soliton. \quad $\Box $

\medskip

A last generalization eliminates the scalar $\lambda $. Recall that given the pair $(\tilde{\nabla }, F)$ as in Definition 1 {\it the exterior covariant derivative of $F$ with respect
to} $\tilde{\nabla }$ is:
$$
(d^{\tilde{\nabla }}F)(X, Y):=(\tilde{\nabla }_XF)Y-(\tilde{\nabla }_YF)X+F(\tilde{T}(X, Y)). \eqno(37)
$$
Since $(24)$ is expressed as:
$$
d^{\tilde{\nabla }}\xi +F+\lambda I=0 \eqno(38)
$$
we introduce:

\medskip

{\bf Definition 4} Let $(M, \tilde{\nabla }, F)$ be given with a symmetric $\tilde{\nabla }$. The vector field $\xi $ will be a $(\tilde{\nabla }, F)$-{\it weak soliton} if:
$$
d^{\tilde{\nabla }}\circ  d^{\tilde{\nabla }}(\xi )+d^{\tilde{\nabla }}F=0. \eqno(39)
$$

\smallskip

{\bf Remark 3} For the Levi-Civita connection $\nabla $:
$$
d^{\nabla }\circ d^{\nabla }(Z)=Riem(\cdot , \cdot )Z, \eqno(40)
$$
while Lemma 2 of \cite[p. 182]{mo:r} gives:
$$
d^{\nabla }Q=-d_2\circ \delta ^{\nabla }Riem. \eqno(41)
$$
with the right hand side:
$$
(d_2\circ \delta ^{\nabla }Riem)(X, Y)=\sum _{k=1}^n(\nabla _{e_k}Riem)(X, Y, e_k) \eqno(42)
$$
for an arbitrary orthonormal field $e=\{e_k\}_{k=1,..., n}$ on $(M, g)$. Hence $\xi $ is a Ricci weak-soliton if and only if:
$$
Riem(\cdot , \cdot)\xi =d_2\circ \delta ^{\nabla }Riem (\cdot , \cdot). \eqno(43)
$$
\quad $\Box $

\bigskip

\noindent Faculty of Mathematics \newline University "Al. I.Cuza" \newline Iasi, 700506 \newline Rom\^{a}nia \newline
E-mail: mcrasm@uaic.ro \\
 http://www.math.uaic.ro/$\sim$mcrasm

\end{document}